\magnification 1200
\font \petit = cmr8 at 8pt
\font \ppetit = cmbx8 at 8pt
\font \Petit = cmitt10 at 8pt
\font \Grand = cmbx12 at 16pt
\font \grand = cmr10 at 12pt

\let \n = \noindent

\let \u = \underbar

\let \ce = \centerline

\let \h = \hbox
\def \trait (#1) (#2) (#3){\vrule width #1pt height #2pt depth #3pt}

\def\fudge{\mathchoice{}{}{\mkern.5mu}{\mkern.8mu}} 
\def\bbc#1#2{{\rm \mkern#2mu\vbar\mkern-#2mu#1}} 
\def\bbb#1{{\rm I\mkern-3.5mu #1}} \def\bba#1#2{{\rm #1\mkern-#2mu\fudge #1}} 
\def\bb#1{{\count4=`#1 \advance\count4by-64 \ifcase\count4\or\bba A{11.5}\or    
\bbb B\or\bbc C{5}\or\bbb D\or\bbb E\or\bbb F \or\bbc G{5}\or\bbb H\or    
\bbb I\or\bbc J{3}\or\bbb K\or\bbb L \or\bbb M\or\bbb N\or\bbc O{5} \or    
\bbb P\or\bbc Q{5}\or\bbb R\or\bbc S{4.2}\or\bba T{10.5}\or\bbc U{5}\or    
\bba V{12}\or\bba W{16.5}\or\bba X{11}\or\bba Y{11.7}\or\bba Z{7.5}\fi}} 
\def \R {{\bb R}}

\def \trait (#1) (#2) (#3){\vrule width #1pt height #2pt depth #3pt}
\def \fin{\hfill
	\trait (0.1) (5) (0)
	\trait (5) (0.1) (0)
	\kern-5pt
	\trait (5) (5) (-4.9)
	\trait (0.1) (5) (0)
\medskip}

\parindent=0,3cm

\ce{\Grand Asymptotic behavior of the unbounded solutions}
\vskip 0,2cm
\ce{\Grand of some boundary layer equations}
\vskip 0,5cm
\ce{by}
\vskip 0,5cm
\ce{\grand Bernard BRIGHI and Jean-David HOERNEL}
\vskip 1,2cm

{\ppetit{Abstract.}}{\petit{ We give an asymptotic equivalent at infinity of the
unbounded solutions of some boundary layer equations arising in fluid mechanics.}}

{\footnote{}{\n{\Petit{Mathematics Subject Classification}} {\petit{(2000):} 
{\petit{ 34B40, 34E99.}}}}}

\vskip 1cm
\n Let us consider the following boundary layer differential equation
$$f'''+ff''-\beta f'^2=0\eqno(1)$$
where $\beta<0$. We are interested in non constant solutions (that we will simply call 
{\sl solutions}) of
$(1)$ defined on some interval
$[t_0,\infty)$ and such that 
$$f'(\infty):=\lim_{t\to\infty}f'(t)=0.\eqno(2)$$ 
Equation $(1)$ can be obtained from similarity boundary layer equations as those introduced
by numerous authors in [1], [2], [11], [12], [13], [14], [17] and
[18], and studied from mathematical point of view in [3], [4], [6], [7], [9], [10] and [15].
In these papers, the corresponding differential equation is considered on $[0,\infty)$ with
the boundary conditions
$f(0)=a$, $f'(0)=1$ and $(2)$, or $f(0)=a$, $f''(0)=-1$ and $(2)$. Here, we will be
concerned by unbounded solutions of these problems, and to be as general as possible we
will  consider all the unbounded solutions of $(1)$-$(2)$ defined on some interval
$[t_0,\infty)$. The restriction to $\beta<0$ is due to the fact that for $\beta\geq 0$ none
of the solutions of $(1)$-$(2)$ are unbounded (see Remark 6 below).

For $\beta=0$, equation $(1)$ reduces to the Blasius equation and a lot of papers have
been published about it. To have a survey, we refer to [16], [5], [8] and the references
therein.

Concerning the existence of unbounded solutions of $(1)$-$(2)$, elementary direct methods
give it for $-2\leq\beta<0$ (see for example [7] and [15]). It seems more difficult to get
such existence results for $\beta<-2$ and the best way to overcome this difficulty should
consist in introducing appropriate blow-up coordinates. Precisely, if $f$ is a solution of
$(1)$ which does not vanish on some interval $I$, we set
$$\forall t\in I,~~~~~s=\int_\tau^tf(\xi)d\xi,~~~~~u(s)={f'(t)\over f(t)^2}~~~~~\hbox{ and
}~~~~~v(s)={f''(t)\over f(t)^3}.$$
Then, we easily get 
$$\left\{\matrix{\dot u=v-2u^2,\hfill\cr\noalign{\vskip1mm}
\dot v=-v+\beta u^2-3uv,}\right.$$
where the dot is for differentiating with respect to the variable $s$.
The plane dynamical system that we obtain has the origin as a saddle-node, and studying
the phase portrait in the neighbourhood of it allows us to underscore the fact that unbounded
positive or negative solutions of $(1)$-$(2)$ have to exist. For details, see [9] or [10].

We now focus our attention on the behavior at infinity of these unbounded solutions. We
start by some elementary and useful lemmas.
\vskip 0,5cm
{\bf Lemma 1.} {\it Let $f$ be a solution of
$(1)$ defined on some interval $J$. If there is $\tau\in J$ such that $f''(\tau)\leq 0$, then
for all $t\in J$ such that $t>\tau$ we have $f''(t)<0$.}
\vskip 0,2cm
{\tt Proof.} This follows immediately from the equality
$(f''e^F)'=\beta f'^2e^F,$
where $F$ denotes any anti-derivative of $f$ on $J$, and from the fact that $f'$ and $f''$ cannot
vanish together without $f$ being constant.
\fin

\vskip 0,5cm
{\bf Lemma 2.} {\it Let $f$ be a solution of
$(1)$-$(2)$ defined on some interval $[t_0,\infty)$. There exists $t_1\geq t_0$ such that
$f''(t)f'(t)<0$ and $f'''(t)f''(t)<0$ for $t\geq t_1$.}
\vskip 0,2cm
{\tt Proof.} By Lemma 1, we know that $f''$ cannot vanish more than once on $[t_0,\infty)$
and thus there exists $t_2\geq t_0$ such that 
$$\forall t>t_2,~~~~~f''(t)f'(t)=-f''(t)\int_t^\infty f''(s)ds<0.$$
Differentiating $(1)$ we get 
$f^{(iv)}+ff'''-(2\beta-1)f'f''=0$
and
$(f'''e^F)'=(2\beta-1) f''f'e^F$
where $F$ denotes any anti-derivative of $f$ on $[t_2,\infty)$. It follows that $f'''$ cannot
vanish more than once on $[t_2,\infty)$ in such a way that $f''(t)\to 0$ as $t\to\infty$ and
there exists
$t_1\geq t_2$ such that 
$$\forall t\geq t_1,~~~~~f'''(t)f''(t)=-f'''(t)\int_t^\infty f'''(s)ds<0.$$
This completes the proof.\fin
\vskip 0,5cm

We now are able to prove our main result.

\vskip 0,5cm
{\bf Theorem 3.} {\it Let $f$ be an unbounded solution of
$(1)$-$(2)$. There exists a constant $c>0$ such that  
$$\vert f(t)\vert\sim c t^{1\over 1-\beta}~~~~~\h{ as }~~~~t\to\infty.\eqno(3)$$}

\n{\tt Proof.} Let $f:[t_0,\infty)\to\R$ be an unbounded solution of
$(1)$-$(2)$. 
\vskip 0,2cm
\n\u{Case 1}. Let us assume first that $f$ is positive at infinity. Thanks to Lemma 2,
there exists
$t_1\geq t_0$ such that 
$$\forall t\geq t_1,~~~~~f(t)>0,~~~f'(t)>0,~~~f''(t)<0~~\h{ and }~~f'''(t)>0.$$ 
Therefore, on $(t_1,\infty)$, we have
$(f'f^{-\beta})'=(ff''-\beta f'^2)f^{-\beta-1}=-f'''f^{-\beta-1}<0$
in such a way that the function $\phi=f'f^{-\beta}$ is decreasing on $[t_1,\infty)$ and 
$$\phi(t)=f'(t)f^{-\beta}(t)\longrightarrow l_0\in[0,\infty)~~~\h{ as
}~~~t\to\infty.\eqno(4)$$ 
Now, multiplying equation $(1)$ by $f^{-\beta-1}$
and integrating between $s\geq t_1$ and $t\geq s$ we easily get 
$$\openup2mm\eqalignno{f^{-\beta-1}(t)f''(t)-f^{-\beta-1}(s)f''(s)&+f^{-\beta}(t)f'(t)
-f^{-\beta}(s)f'(s)\cr
&=-(\beta+1)\int_{s}^tf^{-\beta-2}(r)f'(r)f''(r)dr.&(5)}$$
Since $ff'f''<0$ on $(t_1,\infty)$, the right hand side of $(5)$ has a limit as $t\to\infty$
and thus from $(4)$ we deduce that $f^{-\beta-1}(t)f''(t)$ has a limit $l_1\in[-\infty,0]$
as $t\to\infty$. Suppose now $l_1<0$. 
Then there exist $l_2\in(l_1,0)$ and $t_2>t_1$ such
that $f^{-\beta-1}(t)f''(t)<l_2$ for $t>t_2$. It follows that 
$$\forall t>t_2,~~~~~f''(t)<l_2f^{\beta+1}(t)<{l_2\over\phi(t_1)}f'(t)f(t).$$
Integrating, we get
$$\forall t>t_2,~~~~~f'(t)-f'(t_2)<{l_2\over2\phi(t_1)}(f^2(t)-f^2(t_2))$$
and a contradiction with $(2)$ since the right hand side tends to $-\infty$ as $t\to\infty$.
Consequently $l_1=0$ and coming back to $(5)$ we get 
$$l_0=f^{-\beta-1}(s)f''(s)+f^{-\beta}(s)f'(s)
-(\beta+1)\int_{s}^\infty f^{-\beta-2}(r)f'(r)f''(r)dr,\eqno(6)$$
and this equality holds for all $s\geq t_1$. It remains to show that $l_0>0.$ For that we
have to distinguish between the cases
$\beta\geq -1$ and $\beta<-1$.
\vskip 0,1 cm
Assume first that $\beta\geq -1$. Then $(6)$ implies that 
$$l_0\geq \sup_{s\geq t_1}\{f^{-\beta-1}(s)f''(s)+f^{-\beta}(s)f'(s)\}>0$$
because, on the contrary, we should have
$f''(s)+f(s)f'(s)\leq 0$ for all $s\geq t_1$,
and by integrating
$$\forall s\geq t_1,~~~~~f'(s)+{1\over 2}f^2(s)\leq f'(t_1)+{1\over 2}f^2(t_1)$$
which is absurd since $f(s)\to\infty$ as $s\to\infty$.
\vskip 0,1 cm
Assume now that $\beta<-1$. Since the function $\phi$ is decreasing, we have
$$\int_{s}^\infty f^{-\beta-2}(r)f'(r)f''(r)dr\geq
f^{-\beta-2}(s)f'(s)\int_{s}^\infty f''(r)dr=-f^{-\beta-2}(s)f'(s)^2.$$
We then deduce from $(6)$ that
$$\forall s\geq
t_1,~~~~~l_0\geq f^{-\beta-2}(s)\{f'(s)f^{2}(s)+f''(s)f(s)+(\beta+1)f'(s)^2\}.$$ Looking next
at the polynomial $P_s(X)=f'(s)X^2+f''(s)X+(\beta+1)f'(s)^2$, we easily see that for 
$$X>-{f''(s)\over f'(s)}+\sqrt{-(\beta+1)f'(s)},$$ 
we have $P_s(X)>0$. To conclude,
it is sufficient to remark that there exists $s_0\geq t_1$ such that
$$f(s_0)>-{f''(s_0)\over f'(s_0)}+\sqrt{-(\beta+1)f'(s_0)}.\eqno(7)$$
Indeed, on the contrary we should have 
$f'(s)\to\infty$ as $s\to\infty$ 
and a contradiction.
Therefore $(7)$ holds and we have $l_0\geq f^{-\beta-2}(s_0)P_{s_0}(f(s_0))>0$.
\vskip 0,1 cm
Finally, we have $f'(t)f(t)^{-\beta}\sim l_0$ as $t\to\infty$, and by integrating we obtain
$$f(t)^{-\beta+1}\sim l_0(1-\beta) t~~~\h{ as }~~~t\to\infty$$
and the result in this case.

\vskip 0,2cm
\n\u{Case 2}. Let us assume now that $f$ is negative at infinity. Thanks to Lemma 2,
there exists
$t_1\geq t_0$ such that 
$$\forall t\geq t_1,~~~~~f(t)<0,~~~f'(t)<0,~~~f''(t)>0~~\h{ and }~~f'''(t)<0.$$ 
Then, on $[t_1,\infty)$, we have
$(f'(-f)^{-\beta})'=(-ff''+\beta f'^2)(-f)^{-\beta-1}=f'''(-f)^{-\beta-1}<0$
in such a way that the function $\psi=f'(-f)^{-\beta}$ is decreasing and we have 
$$\psi(t)=f'(t)(-f(t))^{-\beta}\longrightarrow l_0\in[-\infty, 0)~~~\h{ as
}~~~t\to\infty.\eqno(8)$$ 
To conclude as in the first case, it is sufficient to prove that $l_0$ is finite.
Multiplying equation $(1)$ by $(-f)^{-\beta-1}$
and integrating between $s\geq t_1$ and $t\geq s$ we easily get 
$$\openup2mm\eqalignno{f''(t)(-f(t))^{-\beta-1}-f''(s)(-f(s))^{-\beta-1}
&-f'(t)(-f(t))^{-\beta}
+f'(s)(-f(s))^{-\beta}\cr
&=(\beta+1)\int_{s}^t f'(r)f''(r)(-f(r))^{-\beta-2}dr.&(9)}$$
If $\beta\geq-1$, then the right hand side of $(9)$ is non positive, and we see immediatly
that
$l_0$ has to be finite. Let us assume now that
$\beta<-1$ and choose $s$ such that $f'(s)f(s)^{-2}>{1\over 2(\beta+1)}$.
Since the function $\psi$ is decreasing and negative, we get  
$$\int_{s}^t f'(r)f''(r)(-f(r))^{-\beta-2}dr\geq
\psi(t)f(s)^{-2}\int_{s}^t f''(r)dr\geq-{1\over 2(\beta+1)}\psi(t).$$
Then, setting $C(s)=f''(s)(-f(s))^{-\beta-1}-f'(s)(-f(s))^{-\beta}$, we easily deduce 
from $(9)$ that 
$$f''(t)(-f(t))^{-\beta-1}-\psi(t)-C(s)\leq -{1\over 2}\psi(t),$$
which gives $\psi(t)\geq-2C(s)$. Thus $l_0$ is finite.  
\fin
 
\vskip 0,5cm
{\bf Remark 4.} If $\beta=-1$, then $\vert l_0\vert=f''(t_1)+f(t_1)f'(t_1)$ and
$\vert f(t)\vert\sim\sqrt{2\vert l_0\vert t}$ as $t\to\infty$.
\vskip 0,5cm

{\bf Remark 5.} For any $\beta\in\R$ and any $\tau\in\R$, the function
$$t\longmapsto{6\over(2-\beta)(t-\tau)}$$ is a bounded convex solution of $(1)$-$(2)$ on
$[t_0,\infty)$ for all $t_0>\tau$. Bounded concave solutions of $(1)$-$(2)$ exist too (see
[4], [5], [6], [7], [9] and [10]). 
\vskip 0,5cm

{\bf Remark 6.} For $\beta\geq 0$, the solutions of $(1)$-$(2)$ are always bounded.
In fact, suppose that $f:[t_0,\infty)\to\R$ is a solution of $(1)$-$(2)$, then we have 
$f'''\geq-ff''$ in such a way that if $f<0$ at infinity, we deduce from Lemma 2 that there
exists $t_1\geq t_0$ such that necessarily $f''<0$ and $f'''>0$ on
$[t_1,\infty)$. Such a $f$ is bounded. 

If now $f>0$ at infinity and is unbounded, then  $f(t)\to\infty$ as
$t\to\infty$ and there exists $t_1\geq t_0$ such that $f''<0$ and $f>1$ on $[t_1,\infty)$.
Therefore $f'''(t)\geq -f''(t)$ for $t\geq t_1$, and by integrating between $s\geq
t_1$ and $\infty$ we obtain
$-f''(s)\geq f'(s)$. Integrating next between $t_1$ and $t\geq t_1$, we get 
$f'(t_1)-f'(t)\geq f(t)-f(t_1)$ and a contradiction by passing to the limit as $t\to\infty$.

\vskip 0,6cm

\ce{\bf References}
\vskip 0,3cm
{\petit{\n [1] W. H. H. BANKS, Similarity solutions of the boundary layer equations
for a stretching wall. J. de M\'ecan. th\'eor. et appl. {\ppetit 2}, 375-392 (1983). 
\vskip 0,1cm

\n [2] W. H. H. BANKS and M. B. ZATURSKA, Eigensolutions in boundary layer flow
adjacent to a stretching wall. IMA J. Appl. Math. {\ppetit 36}, 263-273 (1986). 
\vskip 0,1cm

\n [3] Z. BELHACHMI, B. BRIGHI, J.M. SAC EP\'EE and K. TAOUS, 
Numerical simulations of free convection about vertical flat plate embedded in 
porous media. Computational Geosciences {\ppetit 7}, 137-166 (2003).
\vskip 0,1cm

\n [4] Z. BELHACHMI, B. BRIGHI and K. TAOUS, Solutions similaires pour
un  probl\`eme de couches 
limites en milieux poreux.
C. R. M\'ecanique {\ppetit 328}, 407-410 (2000).
\vskip 0,1cm

\n [5] Z. BELHACHMI, B. BRIGHI and K. TAOUS, On the concave solutions of the Blasius 
equation.
Acta Math. Univ. Comenianae {\ppetit 69} (2), 199-214 (2000).
\vskip 0,1cm

\n [6] Z. BELHACHMI, B. BRIGHI and K. TAOUS, On a family of differential equations
for boundary layer approximations in porous media.
Euro. J. Appl. Math. {\ppetit 12}, 513-528 (2001).
\vskip 0,1cm

\n [7] B. BRIGHI, On a similarity boundary 
layer  equation. Zeitschrift f\"ur Analysis und ihre Anwendungen, {\ppetit 21} (4),
931-948 (2002).
\vskip 0,1 cm

\n [8] B. BRIGHI, The Crocco change of variable for the Blasius equation. In
preparation.
\vskip 0,1 cm

\n [9] B. BRIGHI and J.-D. HOERNEL , On similarity solutions for boundary layer flows
with prescribed heat flux. To appear in Mathematical Methods in the Applied Sciences.
\vskip 0,1 cm

\n [10] B. BRIGHI and T. SARI, Blowing-up coordinates for a similarity boundary 
layer  equation. To appear in Discrete and Continuous Dynamical Systems (Serie A).
\vskip 0,1 cm

\n [11] M. A. CHAUDARY, J. H. MERKIN and I. POP, Similarity solutions in free
convection boundary layer flows adjacent to vertical permeable surfaces in
porous media. I : Prescribed surface temperature. Eur. J. Mech. B-Fluids {\ppetit 14},
217-237 (1995).
\vskip 0,1cm

\n [12] M. A. CHAUDARY, J. H. MERKIN and I. POP, Similarity solutions in free
convection boundary layer flows adjacent to vertical permeable surfaces in
porous media. II : Prescribed surface heat flux. Heat and Mass Transfer {\ppetit 30},
Springer-Verlag, 341-347 (1995).
\vskip 0,1cm

\n [13] P. CHENG and W. J. MINKOWYCZ, Free Convection About a Vertical Flat Plate
Embedded  in a Porous 
Medium With Application to Heat Transfer From a Dike. J. Geophys. Res. 
{\ppetit 82} (14), 
2040-2044 (1977).  
\vskip 0,1cm

\n [14] H. I. ENE and D. POLI$\check{\h{S}}$EVSKI,  Thermal Flow in Porous Media. D. Reidel 
Publishing Company, Dordrecht, 1987.
\vskip 0,1cm

\n [15] M. GUEDDA, Nonuniqueness of solutions to differential equations for boundary
layer approximations in porous media, C. R. M\'ecanique {\ppetit 330},
279-283 (2002).
\vskip 0,1cm

\n [16] P. HARTMANN,  Ordinary Differential Equations. Wiley, New York, 1964.
\vskip 0,1cm

\n [17] D. B. INGHAM and S. N. BROWN, Flow past a suddenly heated vertical plate in a
porous medium. Proc. R. Soc. Lond. A {\ppetit 403}, 51-80 (1986).
\vskip 0,1cm

\n [18] E. MAGYARI and B. KELLER, Exact solutions for self-similar boundary layer
flows induced by 
permeable stretching walls. Eur. J. Mech. B-Fluids {\ppetit 19}, 109-122 (2000).}}
\vskip 0,9cm

\n Anschrift der Autoren:
\vskip 0,2cm
\n Bernard Brighi

\n Jean-David Hoernel

\n Universit\'e de Haute Alsace

\n Facult\'e des Sciences et 
Techniques

\n 4 rue des fr\`eres Lumi\`ere
 
\n F-68093 Mulhouse Cedex

\n bernard.brighi@uha.fr

\bye